\documentclass[a4paper,12pt]{article}

\usepackage{amsfonts}
\usepackage{amscd,color}
\usepackage{amsmath,amsfonts,amssymb,amscd}
\usepackage{indentfirst,graphicx,epsfig}
\usepackage{graphicx}
\input{epsf}
\usepackage{graphicx}
\usepackage{epstopdf}
\usepackage{caption}

\newtheorem{thm}{Theorem}[section]

\newtheorem{lem}[thm]{Lemma}

\newcommand{\ml}{l\kern-0.55mm\char39\kern-0.3mm}

\def\qed{\hfill \nopagebreak\rule{5pt}{8pt}}

\title{\textbf{Graphs with conflict-free connection \\
number two\footnote{Supported by NSFC No.11531011.}}}
\author{\small \ Hong Chang$^{1}$, \ Trung Duy Doan$^{2,3}$ \thanks{Financial support by the Free State of Saxony (Landesstipendium) is thankfully acknowledged.} , \ Zhong Huang$^{1}$,\\
\small \ Stanislav Jendro{\ml} $^{4}$ \thanks {This work was supported by the Slovak Research and Development Agency under the contract No.\! APVV-15-0116 and by the Slovak VEGA Grant 1/0368/16.}
\ Xueliang Li$^{1}$,
\ Ingo Schiermeyer$^{2}$ \thanks{Part of this research was done while the author was visiting the Center for Combinatorics. Financial support is gratefully acknowledged.}\\[0.2cm]
\small $^{1}$Center for Combinatorics and LPMC \\
\small Nankai University, Tianjin 300071, China\\[0.2cm]
\small $^{2}$Institut f\"ur Diskrete Mathematik und Algebra \\
\small Technische Universit\"at Bergakademie Freiberg \\
\small 09596 Freiberg, Germany\\[0.2cm]
\small $^{3}$ School of applied Mathematics and Informatics \\
\small Hanoi University of Science and Technology, Hanoi, Vietnam \\
\small $^{4}$ Institute of Mathematics, P.~J.~\v Saf\' arik University\\
\small Jesenn\'a 5, 040 01 Ko\v sice, Slovakia\\
{\small Email: changh@mail.nankai.edu.cn, trungdoanduy@gmail.com, stanislav.jendrol@upjs.sk  }\\
{\small  2120150001@mail.nankai.edu.cn, lxl@nankai.edu.cn,  Ingo.Schiermeyer@tu-freiberg.de}\\
}
\date{February 12, 2018}
\begin{document}
\maketitle
\begin{abstract}
An edge-colored graph $G$ is \emph{conflict-free connected} if
any two of its vertices are connected by a path, which contains
a color used on exactly one of its edges. The
\emph{conflict-free connection number} of a connected graph $G$,
denoted by $cfc(G)$, is the smallest number of colors needed in
order to make $G$ conflict-free connected. For a graph $G,$ let $C(G)$ be the subgraph of $G$ induced by its set of cut-edges.
In this paper, we first show
that, if $G$ is a connected non-complete graph $G$ of order $n\geq 9$ with $C(G)$ being a linear forest and with the minimum degree
$\delta(G)\geq \max\{3, \frac{n-4}{5}\}$, then  $cfc(G)=2$. The bound on the minimum degree is best possible.
Next, we prove that, if $G$ is a connected non-complete graph of order $n\geq 33$ with $C(G)$ being a linear forest and with $d(x)+d(y)\geq \frac{2n-9}{5}$
for each pair of two nonadjacent vertices $x, y$ of $V(G)$, then $cfc(G)=2$. Both bounds, on the order $n$ and the degree sum, are tight.
Moreover, we prove several results concerning relations between degree conditions on $G$ and the number of cut edges in $G$.\\[2mm]
\end{abstract}

\section{Introduction}

All graphs in this paper are undirected, finite and simple. We
follow \cite{BM} for graph theoretical notation and terminology not
described here. Let $G$ be a graph. We use $V(G), E(G), n(G), m(G)$,
and $\delta(G)$ to denote the vertex-set, edge-set, number of
vertices, number of edges, and minimum degree of $G$, respectively.
For $v \in V(G)$, let $N(v)$ denote the neighborhood of $v$ in $G$,
$deg(x)$ denote the degree of $v$ in $G$.

Let $G$ be a nontrivial connected graph with an associated
\emph{edge-coloring} $c : E(G)\rightarrow \{1, 2, \ldots, t\}$, $t
\in \mathbb{N}$, where adjacent edges may have the same color. If
adjacent edges of $G$ are assigned different colors by $c$, then $c$
is a \emph{proper (edge-)coloring}. For a graph $G$, the minimum
number of colors needed in a proper coloring of $G$ is referred to
as the \emph{edge-chromatic number} of $G$ and denoted by
$\chi'(G)$. A path of an edge-colored graph $G$ is said to be a
\emph{rainbow path} if no two edges on the path have the same color.
The graph $G$ is called \emph {rainbow connected} if every pair
of distinct vertices of $G$ is connected by a rainbow path in $G$.
An edge-coloring of a connected graph is a \emph{rainbow connection
coloring} if it makes the graph rainbow connected. This concept of
rainbow connection of graphs was introduced by Chartrand et
al.\cite{CJMZ} in 2008. For a connected graph $G$, the \emph{rainbow connection number}
$rc(G)$ of $G$ is defined as the smallest number of colors
that are needed in order to make $G$ rainbow connected. Readers interested
in this topic are referred to \cite{LSS, LS1, LS2} for a survey.

Inspired by the rainbow connection coloring and the proper coloring in graphs,
Andrews et al.\cite{ALLZ} and Borozan et al.\cite{BFGMMMT} independently
introduced the concept of a proper connection coloring. Let $G$ be a nontrivial
connected graph with an edge-coloring. A path in $G$ is called a \emph{proper path}
if no two adjacent edges of the path receive the same color. An edge-coloring $c$
of a connected graph $G$ is a \emph{proper connection coloring} if every pair of
distinct vertices of $G$ is connected by a proper path in $G$. And if $k$ colors
are used, then $c$ is called a \emph{proper connection $k$-coloring}. An edge-colored
graph $G$ is \emph{proper connected} if any two vertices of $G$ are connected by a
proper path. For a connected graph $G$, the minimum number of colors that are
needed in order to make $G$ proper connected is called the \emph{proper connection number}
of $G$, denoted by \emph{$pc(G)$}. Let $G$ be a nontrivial connected graph of order
$n$ and size $m$ (number of edges). Then we have that $1\leq pc(G) \leq \min\{\chi'(G), rc(G)\}\leq m$.
For more details, we refer to \cite{ABBFKS,GLQ,HLQM,LLZ} and a dynamic survey \cite{LM}.

Our research was motivated by the following three results.

\begin{thm}\label{thm1-1}\upshape\cite{BDS}
If $G$ is a 2-connected graph of order $n = n(G)$ and minimum degree
$\delta(G) > max\{2, \frac{n+8}{20}\},$ then $pc(G) \leq 2.$
\end{thm}
\begin{thm}\label{thm1-2}\upshape\cite{BDS}
For every integer $d \geq 3,$ there exists a 2-connected graph of order
$n = 42d$ such that $pc(G) \geq 3.$
\end{thm}
\begin{thm}\label{thm1-3}\upshape\cite{HLQM}
Let $G$ be a connected noncomplete graph of order $n \geq 5.$
If $G \notin \{G_1, G_2\}$ and $\delta(G) \geq \frac{n}{4},$ then $pc(G) = 2,$ where $G_1$ and
$G_2$ are two exceptional graphs on $7$ and $8$ vertices.
\end{thm}

A coloring of the vertices of a hypergraph $H$ is called
\emph{conflicted-free} if each hyperedge $E$ of $H$ has a vertex of
unique color that is not repeated in $E$. The smallest number
of colors required for such a coloring is called the
\emph{conflict-free chromatic number} of $H$. This parameter was
first introduced by Even et al. \cite{ELR} in a geometric setting,
in connection with frequency assignment problems for cellular
networks. One can find many results on the conflict-free coloring, see
\cite{CKP,CT,PT}.

Recently, Czap et al. \cite{CJV} introduced the concept of
a conflict-free connection of graphs. An edge-colored graph $G$ is
called \emph{conflict-free connected} if each pair of distinct
vertices is connected by a path which contains at least one color
used on exactly one of its edges. This path is called a
\emph{conflict-free path}, and this coloring is called a
\emph{conflict-free connection coloring} of $G$. The
\emph{conflict-free connection number} of a connected graph $G$,
denoted by $cfc(G)$, is the smallest number of colors needed to
color the edges of $G$ so that $G$ is conflict-free connected. In
\cite{CJV}, they showed that it is easy to compute the conflict-free
connection number for $2$-connected graphs and very difficult for
other connected graphs, including trees.

This paper is organized as follows. In Section {2}, we list some
fundamental results on the conflict-free connection of graphs.
In Sections $3$ and {4}, we prove our main results.

\section{Preliminaries}

At the very beginning, we state some fundamental results on the
conflict-free connection of graphs, which will be used in the sequel.

\begin{lem}{\upshape \cite{CJV}}\label{lem2-1}
If $P_n$ is a path on $n$ edges, then $cfc(P)=\lceil
\log_2(n+1)\rceil $.
\end{lem}

\vskip0.3cm

Let $C(G)$ be the subgraph of $G$ induced on the set of
cut-edges of $G$. The following lemmas respectively provide
a necessary condition and a sufficient condition for graphs $G$ with $cfc(G)=2$.

Recall that a linear forest is a forest where each of its components is a path.

\begin{lem}{\upshape \cite{CJV}}\label{lem2-2}
If $cfc(G)=2$ for a connected graph $G$, then $C(G)$ is a
linear forest whose each component has at most three edges.
\end{lem}

\begin{lem}{\upshape \cite{CJV}}\label{lem2-3}
If $G$ is a connected graph, and $C(G)$ is a
linear forest in which each component is of order $2$, then $cfc(G)=2$.
\end{lem}

The following lemma, which can be seen as a corollary of Lemma \ref{lem2-3}
for $C(G)$ being empty, is of extra interest. A rigorous proof can be found in \cite{DLL}.

\begin{lem}{\upshape \cite{CJV, DLL}}\label{lem2-4}
If $G$ is a $2$-edge-connected non-complete graph, then $cfc(G)=2$.
\end{lem}

A \emph{block} of a graph $G$ is a maximal connected subgraph of $G$
that has no cut-vertex. If $G$ is connected and has no cut-vertex,
then $G$ is a block. An edge is a block if and only if it is a
cut-edge, this block is called \emph{trivial}. Therefore, any nontrivial
block is $2$-connected.

\begin{lem}{\upshape \cite{CJV}}\label{lem2-5}
Let $G$ be a connected graph. Then from its every nontrivial block
an edge can be chosen so that the set of all such chosen edges forms
a matching.
\end{lem}

Let $C(G)$ be a linear forest consisting of $k$ ($k\geq 0$) components $Q_1,Q_2, \ldots, Q_k$ with $n_i=|V(Q_i)|$ such that $2 \leq n_1\leq n_2 \leq \cdots \leq n_k$. We
now present a stronger result than Lemma \ref{lem2-3}, which will be important to show our main results.

\begin{thm}\label{th2-6}
If $G$ is a connected non-complete graph with $C(G)$ being a linear forest with $2=n_1=n_2=\cdots=n_{k-1}\leq n_k \leq 4$  or $C(G)$ being edgeless, then $cfc(G)=2$.
\end{thm}
\begin{pf}
If $C(G)$ is edgeless then the theorem is true by Lemma \ref{lem2-4}.
If $C(G)$ a linear forest with at least one edge, then $G$ is a non-complete graph and therefore
$cfc(G)\geq 2$. It remains to verify the converse. Note that one can choose from each nontrivial block an edge so that all the chosen edges create a matching set $S$ by
Lemma \ref{lem2-5}. We define an edge-coloring of $G$ as follows. First, we color all edges from $S$ with color $2$, and the edges in $E(G)\setminus \{S\cup Q_k\}$ with color $1$. Next, we only need to color the edges of $Q_k$. If $n_k=2$, then color the unique edge of $Q_k$ with color $1$. If $n_k=3$, then color two edges of $Q_k$ with colors $1$ and $2$. Suppose
$n_k=4$. It follows that $Q_k$ is a path of order $4$, say $w_1w_2w_3w_4$. We color the two edges $w_1w_2$ and $w_3w_4$ with color $1$, and $w_2w_3$ with color $2$. It is easy to check that this coloring is a conflict-free connection coloring of $G$. Thus, we have $cfc(G)\leq 2$, and hence $cfc(G)=2$.\qed
\end{pf}

\vskip0.3cm

\noindent\textbf{Remark 1:} The following example points out that Theorem \ref{th2-6} is optimal in sense of the number of components with more than two vertices of the linear forest $C(G)$ of a graph $G$.

For $t\geq 3$, let $S_n$ be the graph with $n=5t$ vertices, consisting of the path $P_6=v_0v_1v_2v_3v_4v_5$ with complete graphs $K_t$ attached to the vertices $v_i, i\in\{0,1,4,5\}$ and one more $K_t$ sharing the edge $v_2v_3$ with $P_6$. Observe that $\delta(S_n) = t-1 = \frac{n-5}{5}$, and $C(S_n)$ is a linear forest with two components of order $3$, paths $v_0v_1v_2$ and $v_3v_4v_5$.
In any conflict-free connection coloring of $S_n$ with two colors the edges $v_0v_1$ and $v_1v_2$ (resp.  $v_3v_4$ and $v_4v_5$) receive different colors. But then any $v_0$-$v_5$ path has a conflict.
This means that $cfc(S_n) \geq 3$.
\section{Degree conditions and the number of cut-edges}

\vskip0.3cm

\begin{thm}\label{th3-1}
Let $G$ be a connected  graph of order $n\geq k^2, k\geq 3$. If
$\delta(G)\geq{\frac{n-k+1}{k}}$, then $G$ has at most $k-2$ cut edges.
\end{thm}

\begin{pf}
Assume or the sake of contradiction that $G$ has at least $k-1$ cut edges. Let $B$ be a set of $k-1$ cut edges of $G$. Then the graph $G \setminus B$ has exactly $k$ components
$G_1,\dots,G_k$. Consider the following two cases.

Case 1. For every $j \in [k]$ there is a vertex $v_j \in V(G_j)$ such that $N(v_j) \subseteq V(G_j)$. Then every component $G_j$ has at least $\frac{n-k+1}{k} + 1$ vertices and we have
\[n = |V(G)| = \sum_{j=1}^k |V(G_j)| \geq k \cdot (\frac{n-k+1}{k} + 1) = n+1,\] a contradiction.

Case 2. There exists some $i \in [k]$ such that $N(v) \not\subseteq V(G_i)$ for every vertex $v \in V(G_i).$ Then $a = |V(G_i)| \leq k-1$ and every vertex $v \in V(G_i)$ is incident with a cut edge from $B.$ Let $m_i$ denote the degree sum of all the vertices of $V(G_i)$ within $G[V(G_i) \cup B]$. Then we have \[\frac{n-k+1}{k} \cdot a \leq m_i \leq a \cdot (a-1) + k - 1.\]
This, together with the bounds on $a$, provides
\[ 0 \leq a\cdot (a - 1 - \frac{n-k+1}{k}) + k - 1 \leq (k - 1)\cdot (k-2 - \frac{n-k+1}{k}) + k - 1.\]
This leads to $n\leq k^2 - 1,$ a contradiction.
\qed
\end{pf}

\vskip0.3cm
The next theorem shows that the bound on the minimum degree in Theorem 3.1 cannot be lowered.
\vskip0.3cm

\begin{thm}\label{th3-2}
For every $k\geq3$ and $t\geq3$ there exists a connected $n$-vertex graph $H_n$ with $n=k\cdot t$,
$\delta(H_n)=\frac{n-k}{k}$, and $k-1$ cut edges.
\end{thm}

\begin{pf}
The graph $H_n$ consists of a path $P_k$ on k vertices to every vertex of it a complete graph $K_t$ is attached.\qed
\end{pf}

\vskip0.3cm
The following theorem shows that the bound $k^2$ on the number $n$ of vertices in Theorem 3.1 is best possible.


\vskip0.3cm

\begin{thm}\label{th3-3}
For every $k\geq 3$ there exists a graph $R_n$ on $n=k^2 -1$ vertices with $\delta(R_n)=\frac{n-k+1}{k}$ and $k-1$ cut edges.
\end{thm}

\begin{pf}
The graph $R_n$ is a connected graph consisting of a central block $B_0$, isomorphic to the complete graph $K_{k-1}$, $k-1$ blocks $B_1, \dots, B_{k-1}$,
that are complete graphs on $k$ vertices, and a matching $M$ of $k-1$ cut edges. This matches the vertices of $B_0$ with the remaining blocks.\qed
\end{pf}

\begin{thm}\label{th3-4}
Let $G$ be a connected  graph of order
\[n\geq \max\{k^2+k,\frac{\lfloor\frac{k}{2}\rfloor\cdot k(k-2)+k^2-5k+3}{k - 4}\}, k\geq5.\] If
$\deg(x)+\deg(y)\geq\frac{2n-2k+1}{k}$ for any two non-adjacent vertices $x$ and $y$ of $G$, then $G$ has at most $k-2$ cut edges.
\end{thm}

\begin{pf}
Assume for the sake of contradiction that $G$ has at least $k-1$ cut edges. Let $B$ be a set of $k-1$ cut edges of $G$. Then the graph $G \setminus B$ has exactly $k$ components
$G_1,\dots,G_k$. Consider the following two cases.

Case 1. For every $j \in [k]$ there is a vertex $v_j \in V(G_j)$ such that $N(v_j) \subseteq V(G_j)$.

Case 1.1. Let $k$ be even. Then \[n = |V(G)| = \sum_{j=1}^\frac{k}{2} |V(G_j)\cup V(G_{k-j+1})|\geq \frac{k}{2}\cdot (\frac{2n-2k+1}{k}+ 2) = n + \frac{1}{2},\]
a contradiction.

Case 1.2. Let $k$ be odd. Then, w.l.o.g., we can suppose that $|V(G_k)| \geq \frac{n-k+1}{k} + 1$. Therefore,
\[n = |V(G_k)| +\sum_{j=1}^\frac{k-1}{2} |V(G_j)\cup V(G_{k-j})|\geq \frac{n-k+1}{k} + 1+ \frac{k-1}{2}\cdot (\frac{2n-2k+1}{k}+ 2)\] \[ = n + \frac{k+1}{2k},\]
a contradiction.

Case 2. There exists some $i \in [k]$ such that $N(v) \not\subseteq V(G_i)$ for every vertex $v \in V(G_i).$

Case 2.1. There exists only one $i \in [k]$ such that all vertices $v \in V(G_i)$ have $N(v) \not \subseteq V(G_i)$.
Observe that $|V(G_i)|=a \leq k-1$. Notice that every vertex $v \in V(G_i)$ is incident with an edge from $B$, and there is a vertex $y \in V(G_i)$ with $\deg(y) \leq a-1 + \frac{k-1}{a}$.
For any component $G_j, j \neq i \in [k]$, there is
\[|V(G_j)| \geq  \lceil\frac{2n-2k+1}{k}\rceil- \deg (y)  + 1 \geq \lceil\frac{2n-2k+1}{k}\rceil - a + 1 - \frac{k-1}{a} + 1.\]
This means that the number of vertices in $G$ is
\begin{align*}
n&= |V(G)| \geq (k-1)\cdot (\lceil\frac{2n-2k+1}{k}\rceil - a + 1 - \frac{k-1}{a} + 1) + a \\
&\geq (k-1)\cdot (\frac{2n-2k+1}{k} - a + 1 - \frac{k-1}{a} + 1) + a.
\end{align*}
After some manipulations we get
\[n \leq \frac{k(k-1)}{k-2}(a\cdot \frac{k-2}{k-1} + \frac{k-1}{a} -\frac{1}{k}).\]
This, together with the bounds on $a$, provides
\[n \leq \frac{k(k-1)}{k-2}(1\cdot \frac{k-2}{k-1} + \frac{k-1}{1} -\frac{1}{k}).\]
The inequality yields
\[n \leq k^2 +k+\frac{1}{k - 2}.\]
Next we check whether $n=k^2+k$ satisfies the original inequality
\[n = |V(G)| \geq (k-1)\cdot (\lceil\frac{2n-2k+1}{k}\rceil - a + 1 - \frac{k-1}{a} + 1) + a.\]
After some manipulations we get
\[k^2+k\geq k^2+2k-2,\] which is impossible.
Then we have
\[n\leq k^2+k-1,\] a contradiction.

Case 2.2 There exists more than one $i \in [k]$ such that all vertices $v \in V(G_i)$ have $N(v) \not \subseteq V(G_i)$. Assume that there exists a pair of non-adjacent vertices $u, w$ with $u \in V(G_{i_1})$ and $w \in V(G_{i_2}).$ It is possible that $i_1=i_2$. Notice that every vertex in such a component is incident with an edge from $B$, and the two vertices $u$ and $w$ are incident with at most one edge from $B$ in common, then $deg(u)+deg(w)-1\leq k-1$.\\
It implies $n\leq \frac{k^2+2k-1}{2}$, a contradiction. Now we get that every vertex in such components is adjacent to the remaining vertices of such components. Hence all possible configurations have been excluded except for two adjacent singletons $\{u\},\{w\}$ as the only such two components $V_{i_1},V_{i_2}$. As $\deg (u)+\deg (w)-1\leq k-1$, w.l.o.g., we assume that $deg(u)\leq \lfloor\frac{k}{2}\rfloor$.
For any component $G_j, j \neq i_1$ or $ i_2$, then

\[|V(G_j)| \geq \frac{2n-2k+1}{k}- \deg (u)  + 1 \geq \frac{2n-2k+1}{k}  -\lfloor\frac{k}{2}\rfloor + 1.\]
This means that the number of vertices in $G$ is
\[n = |V(G)| \geq (k-2)\cdot (\frac{2n-2k+1}{k} -\lfloor\frac{k}{2}\rfloor + 1) + 2.\]
After some manipulations we get

\[n \leq \frac{\lfloor\frac{k}{2}\rfloor\cdot k(k-2)+k^2-5k+2}{k - 4},\] a contradiction.\qed
\end{pf}

\vskip0.3cm

\noindent\textbf{Remark 2:} Observe that the graph $H_n$ of Theorem 3.2 is a good example showing that the bound on the sum of degrees in Theorem 3.4 is tight.

\vskip0.3cm

The next theorem shows that the bound on $n$ cannot be lower than $k^2+k$.

\vskip0.3cm

\begin{thm}\label{th3-5}
For every $k\geq 5$ there exists a graph $D_n$ on $n=k^2 +k-1$ vertices with $\deg(x)+\deg(y)\geq\frac{2n-2k+1}{k}$ for any two non-adjacent vertices
$x$ and $y$ and having $k-1$ cut edges.
\end{thm}

\begin{pf}
Let $D_n$ be a graph consisting of a vertex $v_0$, $k-1$ blocks $B_1, \dots, B_{k-1}$,
that are complete graphs on $k+2$ vertices, and a set $M$ of $k-1$ cut edges joining the vertex of $v_0$ with the $k-1$ blocks $B_1, \dots, B_{k-1}$. Observe that $D_n$ is a connected graph on $k^2+k-1$ vertices such that $\deg(x)+\deg(y)\geq 2k\geq\frac{2n-2k+1}{k}$ for any two non-adjacent vertices $x$ and $y$.\qed
\end{pf}


\section{Degree conditions for $cfc(G)=2$}

\begin{thm}\label{th4-1}
Let $G$ be a connected non-complete graph of order $n\geq 25$. If
$C(G)$ induces a linear forest and $\delta(G)\geq \frac{n-4}{5}$, then $cfc(G)=2$.
\end{thm}

\begin{pf}
Observe that, by Theorem \ref{th3-1}, the subgraph $C(G)$ of any connected graph $G$ with $\delta(G) \geq \frac{n-4}{5}$ contains at most three cut edges. As $C(G)$ is a linear forest, we conclude that $cfc(G)=2$ by Theorem \ref{th2-6}.
\qed
\end{pf}

\vskip0.3cm

\noindent\textbf{Remark 3:}
The graph $S_n$ defined in the end of Section 2 provides a good example showing the tightness of the minimum degree in Theorem \ref{th4-1}.

\vskip0.3cm

Next, we discuss the minimum degree condition for small graphs to have conflict-free connection number $2$.

\begin{thm}\label{th4-2}
Let $G$ be a connected non-complete graph of order $n$, $9\leq n\leq 24$. If
$C(G)$ induces a linear forest and $\delta(G)\geq \max\{3,\frac{n-4}{5}\}$, then $cfc(G)=2$.
\end{thm}

\begin{pf}
We may assume that $C(G)\neq\emptyset$ by Lemma \ref{lem2-4}. Let $C(G)$ consist of $k$ components $Q_1,Q_2, \ldots, Q_k$ with $n_i=|V(Q_i)|$ such that $2 \leq n_1\leq n_2 \leq \cdots \leq n_k$. We may also assume that $ 3 \leq n_{k-1}\leq n_k\leq 4 $ by Lemma \ref{lem2-2} and Theorem \ref{th2-6}. Then $G\backslash (E(Q_{k-1})\cup E(Q_k))$ has at least five components $C_1,C_2,C_3,C_4,C_5$. Since $\delta(G)\geq 3 $, it follows that $|V(C_i)|>3$ for $1\leq i\leq 5$. Notice that at most two vertices in $C_i$ can be contained in $Q_{k-1}\cup Q_k$, then for each $C_i$ there exists a vertex $u_i$ such that $N(u_i)\subseteq V(C_i)$ for $1\leq i \leq 5$.
Thus, $|V(G)|\geq \sum^{5}_{i=1} |V(C_i)|\geq \sum^{5}_{i=1} (d(u_i)+1) \geq 5(\frac{n-4}{5}+1)=n+1> n$, a contradiction, which completes the proof. \qed
\end{pf}

\vskip0.3cm

\noindent\textbf{Remark 4:} The following examples show that the minimum degree condition in Theorem \ref{th4-2} is best possible. Let $H_i$ be a complete graph of order three for $1\leq i \leq 2$, and take a vertex $v_i$ of $H_i$ for $1 \leq i \leq 2$. Let $H$ be a graph obtained from $H_1,H_2$ by connecting $v_1$ and $v_2$ with a path of order $t$ for $t\geq 5$. Note that $\delta(H)=2$, but $cfc(H)\geq 3 $. Another graph class is given as follows. Let $G_i$ be a complete graph of order $\frac{n}{5}$, and take a vertex $w_i$ of $G_i$ for $1 \leq i \leq 5$. Let $G$ be a graph obtained from $G_1,G_2,G_3,G_4,G_5$ by joining $w_i$ and $w_{i+1}$ with an edge for $1 \leq i \leq 4$. Notice that $\delta(G)=\frac{n-5}{5}$, but $cfc(G)\geq 3$.

\begin{thm}\label{th4-3}
Let $G$ be a connected noncomplete graph of order $n$ with $4 \leq n\leq 8$. If
$C(G)$ induces a linear forest and $\delta(G)\geq 2$, then $cfc(G)=2$.
\end{thm}
\begin{pf}
If $|E(C(G))|\leq 3$, then the proof follows from  Theorem \ref{th2-6}. Otherwise the subgraph $G\setminus E(C(G))$ has at least five components. Since $\delta(G)\geq 2$, at least two components of it have at least three vertices. Thus $|V(G)|\geq 3\times 2+3=9>8$, a contradiction.
\end{pf}

\vskip 0.3cm

\noindent\textbf{Remark 5.} The following example shows that the minimum degree condition in Theorem \ref{th4-3} is best possible. Let $G$ be a path of order $t$ with $t\geq 5$. It is easy to see that $\delta(G)=1$, but $cfc(G)=\lceil \log_2 t\rceil \geq 3 $ by Lemma \ref{lem2-1}.

\vskip 0.3cm

If we do not require that $C(G)$ is a linear forest in above theorems, then we can get the following theorem.

\begin{thm}\label{th4-4}
Let $G$ be a connected non-complete graph of order $n\geq 16$. If
$\delta(G)\geq \frac{n-3}{4}$, then $cfc(G)=2$.
\end{thm}

\begin{pf}
Observe that Theorem \ref{th3-1} shows that $C(G)$ of any connected graph $G$ with $\delta(G)\geq \frac{n-3}{4}$ has at most two edges.
This, when applying Theorem \ref{th2-6}, immediately gives our theorem.
\qed
\end{pf}

\vskip0.3cm

\noindent\textbf{Remark 6:} The following example shows that the minimum degree condition in Theorem \ref{th4-4} is best possible. Let $H_i$ be a complete graph of order $\frac{n}{4}$ for $1\leq i \leq 4$, and take a vertex $v_i$ of $H_i$ for $1 \leq i \leq 4$. Let $H$ be a graph obtained from $H_1,H_2, H_3, H_4$ by adding the edges $v_1v_2, v_1v_3, v_1v_4.$
Note that $\delta(H)= \frac{n-4}{4}$, but $cfc(H)\geq 3 $. On the other hand, the condition $n \geq 16$ in Theorem \ref{th4-4} is also best possible. Let $G_1, G_2, G_3, G_4$ be
complete graphs of order $1, 4, 5, 5$, respectively, and take a vertex $w_i$ of $G_i$ for $1 \leq i \leq 4$. Let $G$ be a graph obtained from  $G_1, G_2, G_3, G_4$ by adding the edges $w_1w_2, w_1w_3, w_1w_4$. Note that $\delta(G)\geq \frac{n-3}{4}$, but $cfc(G) \geq 3$. Also the graph $R_4$ from Theorem \ref{th3-3} shows the sharpness of the bound of $n$.

\begin{thm}\label{th4-5}
Let $G$ be a connected non-complete graph of order $n\geq 33$. If
$C(G)$ is a linear forest, and $\deg(x)+\deg(y)\geq \frac{2n-9}{5}$ for each pair of two non-adjacent vertices $x$ and $y$ of $V(G)$, then $cfc(G)=2$.
\end{thm}

\begin{pf}
From Theorem \ref{th3-4} we deduce that the subgraph $C(G)$ of $G$ has at most three edges. Now  the proof follows from Theorem \ref{th2-6}.
\qed

\end{pf}

\vskip 0.3cm

\noindent\textbf{Remark 7:} An example of the graph $S_n$, introduced in Remark 1, shows that the degree sum condition in Theorem \ref{th4-5} is best possible.
On the other hand, the condition $n\geq 33$ in Theorem \ref{th4-5} is also best possible. Let $G_i$ be a complete graph of order $\frac{n-2}{3}$ for $1 \leq i \leq 3$ and $n \leq 32$, and $G_4 = v_1u_1u_2v_2v_3$ be a path of order 5. Let $G$ be a graph obtained from $G_1, G_2, G_3, G_4$ by identifying a vertex of $G_i$ to the vertex $v_i$ for $1\leq i \leq 3$. Note that the resulting graph $G$ satisfies that $\deg(x)+\deg(y)\geq \frac{2n-9}{5}$ for each pair of two non-adjacent vertices $x$ and $y$ of $V(G)$ and $cfc(G) \geq 3$.


\begin{thebibliography}{1}

\bibitem{ALLZ}
E. Andrews, E. Laforge, C. Lumduanhom, P. Zhang, On proper-path
colorings in graphs, {\it J. Combin. Math. Combin. Comput.}
{\bf 97} (2016), 189--207.

\bibitem{ABBFKS}
S. A. van Aardt, C. Brause, A. P. Burger, M. Frick, A. Kemnitz, and I. Schiermeyer, Proper connection and size of graphs,
{\it Discrete Math.}
{\bf 340 (11)} (2017), 2673--2677.

\bibitem{BM}
J.A. Bondy, U.S.R. Murty, {\it Graph Theory}, GTM $244$, Springer,
$2008$.

\bibitem{BFGMMMT}
V. Borozan, S. Fujita, A. Gerek, C. Magnant, Y. Manoussakis,
L. Montero, Zs. Tuza, Proper connection of graphs,
{\it Discrete Math.} {\bf 312} (2012), 2550--2560.

\bibitem{BDS}
C. Brause, T. Duy Doan, and I. Schiermeyer, Minimum Degree Conditions for the Proper Connection Number of Graphs,
{\it Graphs and Combinatorics} {\bf 33} (2017), 833--843.

\bibitem{CHLMZ}
H. Chang, Z. Huang, X. Li, Y. Mao, H, Zhao, Nordhaus-Gaddum-type theorem for conflict-free
connection number of graphs, arXiv:1705.08316 [math.CO].

\bibitem{CJMZ}
G. Chartrand, G.L. Johns, K.A. McKeon, P. Zhang, Rainbow
connection in graphs, {\it Math. Bohem.} {\bf 133} (2008), 85--98.

\bibitem{CJV}
J. Czap, S. Jendro{\ml}, J. Valiska,  Conflict-free connection of
graphs, {\it Discuss. Math. Graph Theory}, in press.

\bibitem{CKP}
P. Cheilaris, B. Keszegh. D. P\'{a}lv\"{o}igyi,
Unique-maximum and conflict-free coloring for hypergraphs and tree
graphs, {\it SIAM J. Discrete Math.} {\bf 27} (2013), 1775--1787.

\bibitem{CT}
P. Cheilaris, G. T\'{o}th, Graph unique-maximum and
conflict-free colorings, {\it J. Discrete Algorithms} {\bf 9}
(2011), 241--251.

\bibitem{DLL} B. Deng, W. Li, X. Li, Y. Mao, H. Zhao, Conflict-free connection
numbers of line graphs, graphs, {\it Lecture Notes in Computer Science} {\bf No.10627}
(2017), 141--151.

\bibitem{ELR}
G. Even, Z. Lotker, D. Ron, S. Smorodinsky, Conflict-free
coloring of simple geometic regions with applications to frequency
assignment in cellular networks, {\it SIAM J. Comput.} {\bf 33}
(2003), 94--136.

\bibitem{GLQ}
R. Gu, X. Li, Z. Qin, Proper connection number of random graphs,
{\it Theoret. Comput. Sci.} {\bf 609(2)} (2016), 336--343.

\bibitem{HLQM}
F. Huang, X. Li, Z. Qin, C. Magnant, Minimum degree condition for proper connection number 2,
{\it Theoret. Comput. Sci.}, DOI 10.1016/j.tcs.2016.04.042, in press.

\bibitem{LLZ}
E. Laforge, C. Lumduanhom, P. Zhang,  Characterizations of graphs
having large proper connection numbers, {\it Discuss. Math. Graph
Theory} {\bf 36(2)} (2016), 439--453.

\bibitem{LM}
X. Li, C. Magnant, Properly colored notions of connectivity--a
dynamic survey, {\it Theory \& Appl. Graphs} {\bf 0(1)} (2015), Art.
2.

\bibitem{LSS}
X. Li, Y. Shi, Y. Sun, Rainbow connections of graphs: A survey,
{\it Graphs \& Combin.}  {\bf 29} (2013), 1--38.

\bibitem{LS1}
X. Li, Y. Sun,  Rainbow Connections of Graphs, Springer Briefs
in Math., Springer, New York, 2012.

\bibitem{LS2}
X. Li, Y. Sun, An updated survey on rainbow connections of graphs --
a dynamic survey, Theory \& Appl. Graphs {\bf 0(1)}(2017), Art.3.
DOI: 10.20429/tag.2017.000103.

\bibitem{PT}
J. Pach, G. Tardos, Conflict-free colourings of graphs and
hypergraphs, {\it Comb. Probab. Comput.} {\bf 18} (2009), 819--834.


\end{thebibliography}
\end{document}